\documentstyle[12pt]{article}
\newfont{\zahlen}{msbm10 scaled 1200}
\newfont{\Zahlen}{msbm10 scaled 900}
\newcommand{\zz}[1]{\mbox {\zahlen #1}}
\newcommand{\z}[1]{\mbox {\Zahlen #1}}
\newcommand{\pct}{\hfill \rule{6pt}{6pt}}
\newfam\meuffam
\font\twelvmeuf=eufm10 at 12 pt
\font\tenmeuf=eufm10
\font\sevenmeuf=eufm7

\textfont\meuffam=\twelvmeuf
\scriptfont\meuffam=\tenmeuf
\scriptscriptfont\meuffam=\sevenmeuf

\def\germ{\fam\meuffam\tenmeuf}

\def\mg{{\germ m}}

\title{\huge\bf Moduli of simple holomorphic pairs and effective divisors}
\author{\sl Siegmund Kosarew \and \sl Paul Lupascu\thanks{Partially supported
by SNF, nr. 2000-055290.98/1}}
\begin{document}
\newtheorem{thm}{Theorem}[section]
\newtheorem{prop}[thm]{Proposition}
\newtheorem{cor}[thm]{Corollary}
\newtheorem{defin}[thm]{Definition}
\newtheorem{claim}[thm]{Claim}
\newtheorem{lema}[thm]{Lemma} 
\newtheorem{conj}[thm]{Conjecture}
\newtheorem{rem}[thm]{Remark}
\date{}
\maketitle
\begin{quote}
{\sc ABSTRACT:}
In this note we identify two complex structures (one is given by
algebraic geometry, the other by gauge theory) on the set of isomorphism
classes of holomorphic bundles with section on a given compact
complex manifold.\\
In the case of {\it line} bundles, these complex spaces are shown to be
isomorphic to a space of effective divisors on the manifold.

\end{quote}
\section*{Introduction}
Let $(X,{\cal O}_X)$ be a compact complex analytic
space.
We denote by
${\rm Div}^+(X)$ the space  of all Cartier divisors (including the
empty one) on $X$. This set is a Zariski open subspace  in the entire
Douady space
${\cal D}(X)$, parameterizing
{\it all} compact subspaces of $X$. If $X$ is smooth (or more generally
{\it locally factorial}),  then ${\rm Div}^+(X)$ is a union of  connected
components of
${\cal D}(X)$
\cite{fujiki}.\vspace{0.2cm}\\
We consider pairs 
$({\cal E},\phi)$ consisting of an invertible sheaf $\cal E$ over $X$ coupled 
with a holomorphic section 
$\phi$ in $\cal E$, which is locally a non zero divisor.  
Two such  pairs $({\cal E}_i,\phi_i)$ are called 
equivalent if there exists an isomorphism of sheaves
$\Theta:{\cal E}_1\longrightarrow {\cal E}_2$ such that
$\Theta\circ\phi_1=\phi_2$.
The set of equivalence classes is called
the moduli space of simple holomorphic pairs of rank one on $X$.
This moduli space can be given a structure of a complex
analytic space using standard results of deformation theory (see [KO]).\\
There exists a natural bijective
map  from the moduli space of simple holomorphic pairs of rank one into
$ {\rm Div}^+(X)$ sending the equivalence class of a pair $({\cal
E},\phi)$ into the divisor, given by the vanishing locus of the section
$\phi$. In the first part of the paper, we prove that this one-to-one
correspondence is in fact an analytic isomorphism with respect to these
natural structures on the two spaces.\vspace{0.2cm}\\
If $X$ is a 
{\it smooth} compact complex manifold, there is a second
possibility of defining an analytic structure on the moduli space of
simple holomorphic pairs, namely by using gauge-theoretical methods
(compare [OT],[S]). In
the last part we show that, in this case,  these  two structures are
isomorphic.  This is the  "pair"-version of a previous result due to
Miyajima [M] in the case of moduli spaces of simple {\it
bundles}.\newpage   
Throughout this paper we adopt the following notations:
\vspace*{-0.7cm}\begin{tabbing}
mmmm\=m\=mmmmmmmmmmmmmmmmm\kill\\
({\it sets}) \>:\>  the category of sets;\\
({\it an}) \>:\> the category of (not necessarily reduced) complex
spaces;\\ ({\it an}/R)\>:\> the category of relative complex spaces over
a complex space
$R$;\\
({\it germs}) \>:\>  the category of germs of complex spaces;\\
$h_S$ \>:\> the canonical contravariant functor $h_S:{\cal
C}\longrightarrow ({\it sets})$, $h_S={\rm Hom}(\;\cdot\;,S)$ 
\+\+\\associated to an object
$S$  belonging to a category
$\cal C$.
\end{tabbing}
\section{Main result}
The purpose of this part is to prove the
following
\begin{thm}\label{mare}
Let $X$ be a compact complex space with $\Gamma(X,{\cal O}_X)=\zz C$. 
Then there exists a natural complex analytic isomorphism between the
moduli space of simple holomorphic pairs of rank one and  the space ${\rm
{\rm Div}}^+(X)$.
\end{thm}
{\sc Remark.}
Let $(X,{\cal O}_X)$ be a compact complex space with $\Gamma (X,{\cal
O}_X)=\zz C$. Then:\vspace{0.2cm}\\
i.) The Picard functor $\underline{\rm Pic}_X$ is representable (cf.
[Bi1, p.337]).\\ ii.) There exists a Poincar\'e line bundle over
$X\times{\rm Pic}(X)$.\vspace{0.2cm}\\
The existence of a Poincar\'e bundle can be seen as follows (see [Br,
p.55] in the smooth case). The Leray spectral sequence for the projection
morphism
$\pi:X\times S\longrightarrow S$ leads to the exact sequence 
$$
0\to {\rm Pic}(S)\to {\rm Pic}(X\times S)
\to H^0({\cal R}^1\pi_*{\cal O}^*_{X\times S})
\to H^2(\pi_*{\cal O}_{X\times S})\to H^2({\cal
O}^*_{X\times S})\;.
$$
Since $\Gamma (X,{\cal O}_X)=\zz C$, the last morphism
is {\it injective}, hence one gets an exact
sequence
\begin{equation}\label{leray}
0\to {\rm Pic}(S)\to {\rm Pic}(X\times S)
\to H^0({\cal R}^1\pi_*{\cal O}^*_{X\times S})
\cong{\rm Hom}(S,{\rm Pic}(X))\to 0\;.
\end{equation}
For $S={\rm Pic}(X)$, this sequence leads to
a line bundle $\cal P$ on $X\times{\rm Pic}(X)$ such that
${\cal
P}|_{X\times\{[\xi]\}}\cong\xi$ for every $[\xi]\in {\rm {\rm
Pic}}(X)$.\vspace{0.2cm}\pct\\
A complex structure on the space of simple holomorphic 
pairs  is given as follows: 
\\ Let $\cal P$ be a Poincar\'e bundle on $X\times{\rm Pic}(X)$. By [Fl]
there exists a linear fiber space $H$ over ${\rm {\rm Pic}}(X)$,
which represents the functor  $  {\cal H}: (\it an/\mbox{\tiny {\rm
Pic}(X)})\longrightarrow ({\it sets})$ given by
$$\Bigl(S\stackrel{f}{\longrightarrow} {\rm {\rm Pic}}(X)\Bigr)\longmapsto
{\rm Hom}\Bigl({\cal O}_{X\times S}, ({\rm id}_X\times  f)^*{\cal
P}\Bigr)\;.$$ 
(The action of $\cal H$ on morphisms is given by "pull-back".) In
particular, for every complex space
$S$ over
$Pic(X)$, there is a bijection 
\begin{equation}\label{repres}
{\rm Hom}\Bigl({\cal O}_{X\times S}, ({\rm id}_X\times  f)^*{\cal
P}\Bigr)\cong  {\rm Hom}_{Pic(X)}(S,H)\;.
\end{equation}
Let $\tilde O\subset H$ be the subset consisting of zero divisors. The
multiplicative group $\zz C^*$ operates on $H':=H\setminus\tilde O$ such
that the quotient ${\rm P(H')}:=H'/\zz C^*$ becomes an open subset
of a projective fiber space over
${\rm Pic}(X)$.  The fiber over $[\xi]\in {\rm Pic}(X)$ can be identified with
an open subset of 
$PH^0(X,\xi)$. Moreover, ${\rm P(H')}$ coincides set-theoretically with the moduli
space of simple holomorphic pairs, defining a
natural analytic structure on  it.\vspace{0.2cm}\\
In order to prove that ${\rm P(H')}$ and ${\rm Div}^+(X)$ are isomorphic as
complex spaces, it suffices to prove that the associated functors
$h_{{\rm P(H')}}$ and
$h_{{\rm Div}^+(X)}$ are isomorphic. More precisely, we show that both are
isomorphic to the contravariant functor\\
$F:({\it an})\longrightarrow ({\it sets}),$
defined by $S\longmapsto F(S)$, where $F(S)$ denotes the equivalence classes
of pairs $({\cal E},\phi)$, where $\cal E$ is an invertible sheaf on $X\times
S$, and $\phi$ is a holomorphic section in $\cal E$ whose restriction to each
fiber $X\times\{s\}$ is locally a non zero divisor. 
(Two such tuples are called
equivalent if there exists an isomorphism of sheaves $\Theta:{\cal
E}_1\longrightarrow{\cal E}_2$ such that $\Theta\circ\phi_1=\phi_2$.)
\vspace{0.2cm}\\
{\sc Remark.} Note that if $({\cal E},\phi)\sim({\cal E},\phi)$, then
necessarily $\Theta={\rm id}_{\cal E}$. Using this simple observation, it is
easy to see that the functor $F$ is of local nature (it is a sheaf),
i.e.  given any complex space $S$ together with an open covering $S=\cup
S_i$, the following sequence is exact:
 $$F(S)\longrightarrow\prod
F(S_i){\longrightarrow}\prod F(S_i\cap
S_j)\;.$$
\\
Recall that 
a {\it relative Cartier divisor} in $X\times S$ over $S$
is a Cartier divisor $Z\subset X\times S$ which is flat over $S$.
We denote by ${\rm Div}^+_S(X)$ the set of all relative Cartier divisors
(including the empty one) over a fixed
$S$. 
\begin{lema}
Let $S$ be a fixed complex space. The map
$${\cal Z}: F(S)\longrightarrow {\rm Div}^+_S(X)$$ sending the class of
$({\cal E},\phi)$ into the Cartier divisor
given by the vanishing locus of the section $\phi$ is well defined and bijective.
\end{lema}
{\sc Proof.}  It is clear that the vanishing locus 
$Z:=Z(\phi)\subset X\times S$ depends only on the isomorphism
class of
$({\cal E},\phi)$. We need to show  that
$Z$ is flat over $S$. Take $(x,s)\in Z$,
denote $A:={\cal O}_{S,s}$, $B:={\cal O}_{X\times S,(x,s)}$ and let
$\mg\subset A$ be the maximal ideal. Since the section $\phi$ restricted
to each fiber $X\times\{s\}$ is locally a non zero divisor,  we have
${\cal O}_{Z,(x,s)}=B/uB$ for some nonzero divisor
$u$. The flatness of the morphism $A\longrightarrow B/uB$ follows then
from the  general Bourbaki-Grothendieck criterion [Fi, p. 152] since 
$${\rm Tor}_1^A (A/\mg,B/uB)={\rm ker}(B/\mg B\stackrel{\cdot u}{\longrightarrow}
B/\mg B)=\{0\}.$$
We take two simple pairs $({\cal E}_i,\phi_i)$ over
$X\times S$, defining the same divisor $Z:=Z(\phi_1)=Z(\phi_2)$.
The invertible sheaf ${\cal O}_{X\times S}(Z)$ admits a canonical
section $\phi_{can}$,  and both pairs $({\cal E}_i,\phi_i)$ are
equivalent to the pair $({\cal O}_{X\times S}(Z),\phi_{can})$. This
proves the injectivity.\\
Given $Z\in {\rm Div}^+_S(X)$, the associated canonical section
$\phi_{can}$ restricted to each fiber $X\times \{s\}$ is locally a non
zero divisor since $Z$ lies flat over $S$. This gives the surjectivity.
\pct\\
The above lemma shows that there exists
a natural  isomorphism of functors between $F$ and the functor
$G:({\it an})\longrightarrow ({\it sets}),$ $$G(S):={\rm Div}^+_S(X)\subset{\cal
D}(X\times S)\;.$$ It follows from the general result of Douady
\cite{douady} that $G$ is a representable functor and its representation
space is exactly ${\rm Div}^+(X)$. In order to prove theorem \ref{mare} it
suffices to show the following 
\begin{thm}\label{1.4}
The functor $F$ is representable by the space of simple holomorphic pairs
${\rm P(H')}$.
\end{thm}
{\sc Proof.} Since $F$ is a sheaf, it suffices to prove that $F$ and
$h_{{\rm P(H')}}$ are isomorphic as functors defined on the category of 
{germs}  of analytic spaces. The following observation 
shows that $h_{{\rm P(H')}}$ is isomorphic to the sheafified functor
associated to the quotient functor $h_{H'}/h_{\z C^*}$:
\begin{lema} Let $M$ be a complex
analytic space, and let
$G$ be a complex Lie group acting smoothly and freely on  $M$.  Suppose
that the  quotient $M/G$ exists in the category of analytic spaces such that 
the canonical morphism  $M\longrightarrow M/G$ is smooth. Then the 
canonical morphism of functors
\begin{equation}\label{lemaluiK} 
(h_M/h_G)^{\#}\longrightarrow h_{M/G}\;.
\end{equation}
is an isomorphism
(The superscript ${}^{\#}$ denotes here the associated sheafified
functor.)
\end{lema}
{\sc Proof.} Since the projection $M\longrightarrow M/G$ is smooth, one
has an epimorphism of sheaves $h_M\longrightarrow h_{M/G}$, hence the
morphism (\ref{lemaluiK}) is also an epimorphism.\\
Furthermore, since
$G$ acts smoothly and freely on $M$, the morphism $G\times
M\longrightarrow M\times_{M/G}M$,
$(g,m)\longmapsto(m,gm)$ is an {\it isomorphism} by the relativ implicit
function theorem. This shows that (\ref{lemaluiK}) is
also a monomorphism, i.e. an isomorphism.\pct\\
\\
Let $p:H'\longrightarrow {\rm Pic}(X)$ be the natural
morphism and consider  the  corresponding tautological homomorphism
$$u: {\cal O}_{X\times H'}\longrightarrow ({\rm id}_X\times
p)^*{\cal P}$$ given by the bijection (\ref{repres}).\\
The morphism of functors $h_{H'}\longrightarrow F$
is defined by sending  $\phi:S\longrightarrow  H'$
to the isomorphism class of the simple pair
$${\cal O}_{X\times S}\stackrel{\phi^* u}{\longrightarrow}({\rm id}_X\times
p\circ\phi)^*{\cal P}\;.$$
(Note that $\phi^*u|_{X\times\{s\}}=\phi(s)\in H'$ for
every $s\in S$.)\vspace{0.2cm}\\ 
\underbar{\sf Injectivity:} Let $\phi_1,\phi_2$  be two morphisms
from
$S$ to $ H'$ such that the associated simple pairs are
isomorphic.  It follows in particularly that the two sheaves $ ({\rm id}_X\times
p\circ\phi_i)^*{\cal P}$ are isomorphic via some $\Theta$. The 
sequence (\ref{leray}) implies
$p\circ\phi_1=p\circ\phi_2$. The isomorphism $\Theta$ becomes an
automorphism, and is given by multiplication with some element
$a_{\Theta}\in H^0(X\times S,{\cal O}^*_{X\times S})
\cong H^0(S,{\cal O}^*_{S})$ (since $\Gamma(X,{\cal O}_X)=\zz C$). It follows
that
$\phi_1,\phi_2$ are conjugate under the action of
$a_{\Theta}:S\longrightarrow
\zz C^*$, i.e. the morphism 
$$(h_{H'}/h_{\z C^*})^{\#}\longrightarrow F\;$$
is injective.\vspace{0.2cm}\\
\underbar{\sf Surjectivity:} Consider a germ $(S,0)$ of complex space and
a simple pair 
$({\cal O}_{X\times S}\stackrel{\alpha}{\longrightarrow}{\cal E})$.
The corresponding morphism $f\in {\rm Hom}_{Pic(X)} (S,H')$ has the
property $${\cal E}\cong ({\rm id}_X\times p\circ f)^*{\cal
P}\otimes {\cal L}$$ for some $[{\cal L}]\in {\rm Pic}(S)$.
However, since we are working on the germ $S$, we may
assume that ${\cal L}$ is trivial. In this way, we obtain a simple pair 
$({\cal O}_{X\times S}\stackrel{\alpha'}{\longrightarrow}{({\rm id}_X\times
p\circ f)^*{\cal P}})$, which by (\ref{repres}) leads to a morphism 
from $S$ to $H'$. This proves the
surjectivity, and completes the proof of theorem \ref{1.4}.
\pct\\
\begin{cor}\label{fixedclass} Let $(X,{\cal O}_X)$  be a compact, reduced, connected and locally 
irreducible analytic space, and let $m\in H^{2}(X,\zz Z)$ be a fixed
cohomology class.  Consider the closed subspace of ${\rm Div}^+(X)$ given by 
$${\rm Dou}(m):=\Bigl\{\;Z\in {\rm Div}^+(X)\;\Big|\; c_1({\cal
O}_X(Z))=m\;\Bigr\}\;.$$Then there exists an isomorphism of complex
spaces
 $${\rm Dou}(m)\;\cong\Bigl\{\; ({\cal
E},\phi)\;\Big|\;
{\cal E} \mbox{ invertible sheaf
 on $X$, $c_1({\cal E})=m$} ,\;
0\not=\phi\in H^0(X,{\cal E})
\Bigr\}\Big/\sim\;.$$
\end{cor}
{\sc Proof.}  One has a commutative diagram
$$\begin{array}{rcccl}
{\rm P(H')} && \stackrel{\cong}{\longrightarrow} & &{\rm Div}^+(X)\\
&\searrow &  & \swarrow&\\
  & & {\rm Pic}(X) &  &
\end{array}$$
where the vertical arrows are the natural projective morphisms.  The
assertion follows by taking the analytic pull-back of the component 
${\rm Pic}^m(X)\subset {\rm Pic}(X)$ via these two maps.\pct\\
\\
{\sc Remark.} If $X$ is a smooth manifold, then ${\rm Dou}(m)$ is a union of connected components of ${\rm Div}^+(X)$. Moreover, if
$X$ admits a K\"ahler metric, the spaces ${\rm Dou}(m)$ are always
compact. This follows from Bishop's compactness theorem, since all
divisors in
${\rm Dou}(m)$ have the same volume (with respect to any K\"ahler metric). \\
This property fails in the case of manifolds which do not allow
K\"ahler metrics, since a non-K\"ahlerian manifold may
have (nonempty) effective divisors which are homologically trivial.
(Take for instance a (elliptic) surface $X$ with $H^2(X,\zz Z)=0$.)
\section{Gauge-theoretical point of view}
When $X$ is a smooth compact, connected complex manifold it is
possible to construct a "gauge-theoretical" moduli space of simple
holomorphic pairs (of any rank) on $X$ (compare [OT],[S]).\\
\\
 Let $E$ be a
fixed ${\cal C}^{\infty}$ complex vector bundle of rank $r$ on
$X$. We recall the following basic facts from complex differential
geometry:
\begin{defin}{\rm 
A {\it semiconnection} (of type (0,1)) in $E$ is a differential
operator
$\bar\delta: A^{0}(E)\longrightarrow A^{0,1}(E)$ satisfying the Leibniz
rule
$$\bar\delta(f\cdot s)=\bar\partial(f)\otimes
s+f\cdot\bar\delta(s)\qquad\forall\; f\in {\cal C}^{\infty}(X,\zz C), s\in
A^0(E)\;.$$
}
\end{defin} 
The space of all semiconnections in $E$ will be denoted by $\bar{\cal
D}(E)$; it is an affine space over $A^{0,1}({\rm End}E)$. 
Every $\bar\delta\in\bar{\cal
D}(E)$ admits a natural extension 
$\bar\delta: A^{p,q}(E)\longrightarrow A^{p,q+1}(E)$ such that 
$$\bar\delta(\alpha\otimes s)=\bar\partial(\alpha)\otimes s+
(-1)^{p+q}\alpha\wedge\bar\delta(s)\qquad\forall\;\alpha\in
A^{p,q}(X),s\in A^0(E)\;.$$
Moreover, $\bar\delta$ induces $\bar{D}:A^{p,q}({\rm
End}E)\longrightarrow A^{p,q+1}({\rm End}E)$ by
$$\bar{D}(\alpha)=[\bar{\delta},\alpha]:=\bar{\delta}\circ\alpha+(-1)^{p+q+1}
\alpha\circ\bar{\delta}\;.$$
Every {\it holomorphic} bundle
$\cal E$ over
$X$ of differentiable type $E$ induces a canonical
semiconnection
$\bar\delta:=\bar\partial_{\cal E}$ on $E$ such that
$\bar\delta^2:A^{0}(E)\longrightarrow A^{0,2}(E)$ vanishes
identically. Conversely, by [AHS] Theorem (5.1), every
semiconnection
$\bar\delta$ with $\bar\delta^2=0$ defines a unique holomorphic bundle
$\cal E$, differentiably equivalent to $E$ such that $\bar\partial_{\cal
E}=\bar\delta$. \\
There exists a natural
right action of the  gauge group
$GL(E)\subset A^0({\rm End}E)$ of differentiable automorphisms of $E$ on
the space
$\bar{\cal D}(E)\times A^0(E)$ by $$(\bar\delta,\phi)\cdot
g:=(g^{-1}\circ\bar\delta
\circ g,g^{-1}\phi)\;.$$
Denote by $\bar{\cal S}(E)$ the set of points with trivial isotropy
group. After suitable $L^2_k$-Sobolev completions, the space $\bar{\cal
B}^{s.p.}(E):=
\bar{\cal S}(E)/GL(E)$ becomes a complex analytic Hilbert manifold, and 
$\bar{\cal S}(E)\longrightarrow \bar{\cal B}^{s.p.}(E)$ a complex analytic
$GL(E)$-Hilbert principal bundle. The map
$$\Upsilon:\bar{\cal D}(E)\times A^0(E)\longrightarrow
A^{0,2}({\rm End} E)\times A^{0,1}(E)$$ given by 
$$\Upsilon(\bar\delta,\phi)=(\bar\delta^2,\bar\delta\phi)\;$$ 
is $GL(E)$-equivariant, hence it induces a section
$\hat\Upsilon$ in the associated Hilbert vector bundle 
$$\bar{\cal S}(E)\times_{GL(E)}\Bigl[A^{0,2}({\rm End}E )\oplus
A^{0,1}(E)\Bigr]$$ over 
$\bar{\cal B}^{s.p.}(E)$. This section becomes analytic for apropriate
Sobolev completions.
\begin{defin}{\rm The
{\it gauge-theoretical} moduli space
${\cal M}^{s.p.}(E)$ of simple holomorphic pairs of type $E$ is the 
complex analytic space given by the vanishing locus of the
section $\hat\Upsilon$.}
\end{defin}
Set-theoretically,  ${\cal M}^{s.p.}(E)$ can be identified with the set of
isomorphism classes of pairs
$ ({\cal E},\phi)$, where 
${\cal E}$ is a holomorphic bundle of type $E$, and $\phi$ is a 
holomorphic section in $\cal E$, such that the associated evaluation map 
$$ev(\phi): H^0(X,{\cal E}nd({\cal E}))\longrightarrow H^0(X,{\cal E})$$
is injective. This is equivalent to the fact that the only automorphism
of $ ({\cal E},\phi)$ is the identity.
If  $\cal E$ is a
{\it simple bundle} (this happens always if $r=1$) and
$\phi\in H^0(X,\cal E)$,  then $({\cal E},\phi)$ is a {\it simple pair}
iff
$\phi$ is nontrivial. 
\begin{defin}
{\rm Fix
$(\bar{\delta}_0,\phi_0)\in\Upsilon^{-1}(0)$.
A {\it gauge-theoretical family of deformations} of
$(\bar\partial_0,\phi_0)$ parametrized by a germ $(T,0)$ is a complex
analytic map
$$\omega=(\omega_1,\omega_2):(T,0)\longrightarrow (A^{0,1}({\rm End}E),0)
\times(A^0(E),\phi_0)$$
such that the image of the map $(\bar\delta_0+\omega_1,\omega_2)$ is
contained in $\Upsilon^{-1}(0)$, and\\
$(\bar\delta_0+\omega_1,\omega_2):(T,0)\longrightarrow
\Upsilon^{-1}(0)$ is also holomorphic.\\ Two families
$\omega$ and $\omega'$ over $(T,0)$ are called {\it equivalent} if there
exists a complex analytic map
$$g:(T,0)\longrightarrow ({\rm GL}(E),{\rm id}_E)$$
such that $\omega'=\omega\cdot g$.}
\end{defin}
Note that, given such a deformation $\omega=(\omega_1,\omega_2)$, the
family
$\omega_1$ induces uniquely a section in the sheaf $\Bigl({\cal
A}^{0,1}({\rm End}E)\times{\cal A}^{0,0}(E)\Bigr)\hat\otimes{\cal O}_T$,
and conversely. In particular, if
$(T,0)$ is an {\it artinian} germ, then $\omega_1$ induces a morphism
of sheaves 
$${\cal A}^{0,i}(E)_T:={\cal A}^{0,i}(E)\otimes_{\z C}{\cal
O}_T\longrightarrow {\cal A}^{0,i+1}(E)_T$$
We denote by $h_{gauge}:({\it germs})\longrightarrow({\it sets})$ the 
functor which sends a germ $(T,0)$ into the set of
equivalence classes of gauge-theoretical families of deformations over
$(T,0)$.
\begin{thm} The functor $h_{gauge}$ has a semi-universal deformation.
\end{thm}
{\sc Proof.} Fix $\bar\delta_0\in\bar{\cal D}(E)$ and consider the orbit
map
$\beta:GL(E)\longrightarrow\bar{\cal D}(E)\times A^0(E)$ given by
$$\beta(g):=(\bar{\delta}_0,0)\cdot
g=(\bar{\delta}_0+g^{-1}\bar{D}_0(g),0)\;.$$ By [BK] Theorem (12.13) and
[K] Theorem (1.1), the existence of a semi-universal deformation follows
if:
\begin{itemize}
\item[{$i.)$}] The derivative of $\Upsilon$ at $(\bar{\delta}_0,0)$ and
the derivative of $\beta$ at ${\rm id}_E$ are direct linear maps (for
apropriate Sobolev completions);
\item[{$ii.)$}] The quotient ${\rm
ker}(d\Upsilon_{(\bar{\delta}_0,0)})\Big/ {\rm im} (d\beta_{{\rm
id}_E})$ is finite dimensional.
\end{itemize}
One has
$${\rm ker}(d\Upsilon_{(\bar{\delta}_0,0)})=\Bigl\{\;(\alpha,\phi)\in
A^{0,1} (EndE)\times
A^{0}(E)\;\Big|\;\bar{D}_0(\alpha)=0,\bar{\delta}_0\phi=0\;\Bigr\}$$ and
$d\beta_{{\rm id}_E}:A^0(EndE)\longrightarrow A^{0,1} (EndE)\times
A^{0}(E)$ is given by $u\longmapsto (\bar{D}_0(u),0)$. Therefore $i.)$
and
$ii.)$ follow from standard Hodge theory.\pct\\
\\
{\sc Remark.}
The existence of a moduli space of simple holomorphic pairs for
arbitrary rank can be also deduced from [KO] Theorem (2.1), (2.2) and
the proof of Theorem (6.4) of loc.cit.\\ Indeed, one has local
semi-universal deformations: Fix $({\cal E}_0,\phi_0)$ and let ${\cal
E}\longrightarrow X\times (R,0)$ be a semi-universal family of vector
bundles with
${\cal E}|_{X\times\{0\}}\cong {\cal E}_0$. Similarly
as in rank one, the functor ${\cal H}:({\it
an}/R)\longrightarrow({\it sets})$ given by
$$\Bigl(S\stackrel{f}{\longrightarrow} R\Bigr)\longmapsto
{\rm Hom}\Bigl({\cal O}_{X\times S}, ({\rm id}_X\times  f)^*{\cal
E}\Bigr)\;$$ is representable [Fl] by a linear fibre space $p:\tilde
R\longrightarrow R$. Moreover, there exists a tautological section
$u:{\cal O}_{X\times\tilde R}\longrightarrow ({\rm id}_X\times p)^*{\cal
E}$ arising from the representability of $\cal H$.
The pair $(({\rm id}_X\times p)^*{\cal E},u^*\phi_0)$ is a versal
deformation of $({\cal E}_0,\phi_0)$. By [Bi2] there exists
then also a semi-universal deformation. \\
Moreover, it can be shown (as in loc.cit. for the case of
simple sheaves), that the isomorphy locus of two families
of simple holomorphic pairs over $S$ is a {\it locally
closed} analytic subset of $S$ and as a consequence, this
moduli space exists by [KO] (2.1).\vspace{0.2cm}\\ The aim of the
remaining part is to prove the following 
\begin{thm} The {\it
gauge-theoretical} moduli  space of simple holomorphic pairs of type
$E$ is {analytically isomorphic}  to the {\it 
complex-theoretical} moduli space of simple holomorphic pairs of type
$E$. 
\end{thm}
{\sc Proof.} The arguments we use are inspired from [M],
where a similar problem is treated (bundles without section).
It suffices to show that the associated deformation functors
$h_{gauge}$ resp. $h_{an}$ are isomorphic over {\it artinian} bases.
\\ Note first, that there exists a well defined morphism of functors from 
$h_{an}$ to $h_{gauge}$. Moreover, this morphism is injective since
$$({\cal E}_1,\phi_1)\sim({\cal E}_2,\phi_2)\Longleftrightarrow
(\bar\partial_{{\cal E}_1},\phi_1)\sim(\bar\partial_{{\cal
E}_2},\phi_2)\;.$$
In order to prove the surjectivity, we need to show that every
gauge-theoretical family of simple holomorphic pairs determines a
complex-theoretical family of simple holomorphic pairs. 
We will prove this by induction
on the length of the artinian base.\\
For $n=0$, this follows from the fact that every integrable
semi-connection $\bar\delta_0$ determines a holomorphic bundle ${\cal
E}_0$ of type $E$
such that one has an exact sequence of sheaves
$$\begin{array}{cccccccccc}
0 &\to & {\cal E}_0 & \to & {\cal A}^{0,0}(E)
& \stackrel{\bar{\partial}_{{\cal E}_0}}{\to} & {\cal A}^{0,1}(E) &
\stackrel{\bar{\partial}_{{\cal E}_0}}{\to} & {\cal
A}^{0,2}(E)&\to\dots\\
\end{array}$$
For the induction step, let $(T,0)$ be a small extension of an
infinitesimal neighbourhood $(T',0)$ such that 
${\rm ker}({\cal O}_T\longrightarrow {\cal O}_{T'})=\zz C$, and
let $\omega=(\omega_1,\omega_2)$ be a gauge-theoretical family of simple
holomorphic pairs
parametrized by $(T,0)$.
By the induction assumption, we can find a holomorphic vector bundle
${\cal E}'$ over $X\times(T',0)$ which is induced by
$\omega_1|_{T'}$. Then we have the following exact sequence of sheaves
$$
\begin{array}{cccccccccc}
 & & 0 & & 0 & & 0 & &0\\
& &\downarrow & & \downarrow & & \downarrow & &\downarrow\\
0 &\to & {\cal E}_0 & \to & {\cal A}^{0,0}(E)
& \to & {\cal A}^{0,1}(E) & \to & {\cal A}^{0,2}(E)&\to\dots\\
& &\downarrow & & \downarrow & & \downarrow & &\downarrow\\
0 &\to & {\cal E} & \to & {\cal
A}^{0,0}(E)_T & \to & {\cal A}^{0,1}(E)_T &
\to & {\cal A}^{0,2}(E)_T &\to\cdots\\
& &\downarrow & & \downarrow & & \downarrow & &\downarrow\\

0 &\to & {\cal E}' & \to & {\cal
A}^{0,0}(E)_{T'} & \to & {\cal
A}^{0,1}(E)_{T'} &
\to & {\cal A}^{0,2}(E)_{T'}&\to\cdots\\
& &\downarrow & & \downarrow & & \downarrow & &\downarrow\\
 & & 0 & & 0 & & 0 & &0
\end{array}
$$
In particular, $\cal E$ is locally free, hence it defines a holomorphic
vector bundle over $X\times(T,0)$ whose restriction to $X\times(T',0)$
gives ${\cal E}'$.
This vector bundle ${\cal E}$ together with the family of sections
$\omega_2$ gives rise to a complex theoretical family of simple pairs
which induces the gauge-theoretical family $\omega$.\\
This completes the proof of the theorem.\pct

{\bf Key words:} moduli spaces, simple holomorphic pairs,
divisors, semiconnections, holomorphic structures\\
{\bf AMS-classification:} 14 C 20, 14 J 15, 32 G 05, 32 G 13, 53 C 05\\
\\
{\bf Authors addresses:}\vspace*{0.3cm}\\
{\it Siegmund Kosarew}, Institut Fourier, Universit\'e de Grenoble 1,
F-38402 Saint Martin d'H\`eres, {\sc France}\\
{\sl Siegmund.Kosarew@ujf-grenoble.fr}\\
\\
{\it Lupascu Paul}, Mathematisches Institut, Universit\"at Z\"urich,
Winterthurerstr. 190, 8057-Z\"urich, {\sc Switzerland}\\
{\sl lupascu@math.unizh.ch}

\end{document}